\documentstyle[amssymb,amsfonts]{bcp98e}

\newcounter{defn}
\newtheorem{dfn}[defn]{Definition}
\newcounter{thm}
\newtheorem{theorem}[thm]{Theorem}
\newtheorem{conjecture}{Conjecture}
\newtheorem{proposition}[thm]{Proposition}
\newtheorem{lemma}[thm]{Lemma}

\newcommand{\Z}{{\mathbb{Z}}}
\newcommand{\R}{{\mathbb{R}}}

\begin{document}

\mathclass{Primary 17B70; Secondary 81T70, 58A50, 58D29.}

\abbrevauthors{Olga Kravchenko}
\abbrevtitle{Deformations of Batalin--Vilkovisky algebras. }

\title{Deformations of Batalin--Vilkovisky algebras.}

\author{Olga\  Kravchenko}
\address{Institut de Recherche Math\'ematique Avanc\'ee (UMR 7501) \\
          CNRS et Universit\'e Louis Pasteur \\
          7 rue Ren\'e-Descartes   \\
          67084 Strasbourg, France \\
 email ok@alum.mit.edu}

\vfootnote{}{The research for this paper was supported by 
the European Postdoctoral Institut (IPDE).}
\maketitlebcp

\vspace{1mm}

\begin{flushright} {\it To the memory of Stanis\l aw Zakrzewski}
\end{flushright}

\vspace{5mm}

\abstract{We show that a graded commutative algebra $A$ 
with any square zero odd differential operator 
is a natural generalization of a
Batalin--Vilkovisky algebra.
While such an operator of order $2$ defines a Gerstenhaber (Lie) 
algebra structure on $A$,
an operator of an 
order higher than $2$ 
(Koszul--Akman definition) leads to the structure 
of a strongly homotopy Lie  algebra (L$_\infty$--algebra) on $A$. 
This allows us  to give a definition of a 
Batalin--Vilkovisky algebra 
up to homotopy. We also make a
conjecture which is a generalization of the formality 
theorem of Kontsevich to the Batalin--Vilkovisky algebra level. 
}

\section{Introduction.} 

Batalin--Vilkovisky algebras are 
graded commutative algebras with 
an extra structure  given by a second order
 differential operator of square $0$. 
The simplest example  is 
the algebra of polyvector fields on a vector space $\R^n.$
There is a second order square zero
differential operator on this algebra, obtained 
as an operator dual to the de Rham 
differential on the algebra of differential forms \cite{W}.
Namely, if one chooses a volume form, 
one can pair differential forms
to polyvector fields. This pairing lifts 
 the de Rham differential to polyvector 
fields and gives a second order square $0$
operator.

In this article, we consider the following generalization
of the Batalin--Vilkovisky structure:
we do not require that the operator be  of the second order.
The condition that this operator be a differential (of square $0$) 
leads to the structure of $L_\infty$ algebra  \cite{HS,GK,LS} 
(also  called a Lie algebra 
up to homotopy or strong homotopy  Lie  algebra).

The notion of an algebra up to homotopy is a very useful 
tool in proving certain deep  theorems (like 
formality theorem of Kontsevich \cite{K}).

The most important 
property of algebras up to homotopy is 
that all the higher homotopies vanish on their
cohomology groups.
Namely, let $A$ be a $\mathcal P$  algebra up to homotopy,
with a differential $d$; then the space of its cohomology 
$H (A,d)$ is a $\mathcal P$  algebra, where $\mathcal P$  means either
Lie, or  associative, or commutative, or Poisson, or Gerstenhaber, etc.

We propose a definition of a {\it   commutative  strong homotopy
Batalin--Vilkovisky algebra.} Its noncommutative 
version leads to a generalized formality conjecture.

\section{Batalin--Vilkovisky algebras (BV--algebras).}
 We work in the category of $\Z$--graded algebras: 
$A = \oplus A_i.$ 
We denote the degree of 
a homogeneous element $a$ by $|{a}|.$

\begin{dfn}
\label{degree} A map $D\!: A \to A$
is of {\it degree} $|{D}|$ if 
$D\!: A_l \to A_{l+|{D}|}$ for each $l$. 
 
The degree of an 
element $a_1 \otimes \cdots \otimes a_k \in A^{\otimes k}$
is a sum of degrees $\sum_{j=1}^{k}
|{a_j}|.$
\end{dfn} 

Let $\mu\!: A \otimes A \to A$ be a  product on $A$ 
(a priori noncommutative
non-associative).
Following Akman \cite{A}, from any map $ D\!: A \to A$ we can 
inductively define the following 
linear maps $ F_D^{k}\!:  A^{\otimes k} \to A:$
 \begin{eqnarray}
F_D^1 (a)& =&D a,  { } \nonumber  \\
F_D^2 (a_1, a_2)& = & 
D \mu (a_1,a_2) - \mu (Da_1,a_2) -
(-1)^{|{a}_1||{D}|} \mu (a_1, D a_2),  \nonumber\\
 & { \cdots \cdots }  & { } \nonumber  \\
F_D^{n+1} (a_1,..., a_n,a_{n+1}) &  =
&  F_D^n (a_1,...,\mu (a_n, a_{n+1}))  \nonumber \\
\label{eq:FD}          { }   & - & 
 \mu(F_D^n \bigl(a_1,...,a_{n-1},a_n),  a_{n+1})\\
   &  -&(-1)^{|{a}_n|(|{a}_1| + ... +|{a}_{n-1}|
+ |{D}|)}
\mu (a_n, F_D^n (a_1,...,a_{n-1}, a_{n+1})\bigr), \nonumber 
\end{eqnarray}

\begin{dfn} (Akman)\label{Akman}
A linear map $ D\!: A \to A$ is a 
 differential operator of order not higher than  $ k$ if
$F_D^{k+1} \equiv 0.$ 
\end{dfn}

\begin{dfn} A Batalin--Vilkovisky algebra 
(BV--algebra for short) is the following data
 $(A,\delta)\!:$ an  associative 
$\Z$--graded commutative algebra $A$,  
and an operator
 $\delta$ of  order $2$, of degree $-1$, and  of square $0$. 
\end{dfn}

\begin{dfn}
\label{def:G} A Gerstenhaber  algebra is 
a graded space $A = \sum_i A_i$ with 
\begin{itemize}
\item
an associative graded
 commutative product of degree $1, \ \mu: A_i \otimes A_j \to A_{i+j +1}, $

$\mu(a \otimes b)=a \cdot b;$
\item a graded Lie bracket
of degree $0, \ l: A_i \wedge A_j \to A_{i+j}, \  l(a \otimes b) = [a,b],$ 
\item such that the Lie adjoint action 
is an odd derivation with respect to the product:
\[
[a,b \cdot c] = [a,b] \cdot c + (-1)^{|b| \ |c|} [a,c] \cdot b. 
\]
\end{itemize} 
\end{dfn}

\begin{lemma}
Any BV--algebra $(A, \delta)$ is a Gerstenhaber algebra with  
the Lie bracket  given by $F_\delta^2$ up to a sign:
\begin{equation}
\label{eq:BV}
[a_1, a_2] = (-1)^{|{a}_1|} F_\delta^2 (a_1, a_2)= 
(-1)^{|{a}_1|}\bigl( \delta \mu (a_1,a_2) - \mu (\delta a_1,a_2) -
(-1)^{|{a}_1|} \mu (a_1, \delta a_2) \bigr),
\end{equation}
for $a_1, a_2 \in A.$
\end{lemma}
 A Gerstenhaber algebra which is also a BV--algebra is called ``exact''
\cite{KS}, since the bracket then is given by a $\delta$--coboundary.

 \begin{remark}
\rm  In the language of operads one can give  another characterization of a
Gerstenhaber algebra. A Gerstenhaber algebra is an algebra over the braid
operad \cite{G}.  Then BV--algebras are algebras over the {\it cyclic} braid
operad \cite{GK}.  In other words a Gerstenhaber algebra structure 
comes from a BV--operator if the corresponding operad is  cyclic. 
 \end{remark}

\section{$L_\infty$--algebras.}

 The brackets defined by the 
recursive formulas (\ref{eq:FD}) have interesting relations.
We need the  notion of an $L_\infty$--algebra to describe them.

 We view an  $L_\infty$--algebra  structure 
as  a codifferential on the exterior 
coalgebra of a vector space
\cite{LM,P}. . This is a 
generalization of the point of 
view on graded Lie algebras
taken in \cite{R}.

Let $V$ be a graded vector space. Define the 
exterior coalgebra structure  on $\Lambda V$ by giving 
the coproduct  on the exterior algebra
$\Delta : \Lambda V \to \Lambda V \otimes \Lambda V\!:$ 
\begin{eqnarray}
\Delta  v & =& 0  
\label{eq:coproduct}\\
\Delta (v_1 \wedge \cdots \wedge v_n)& = &\sum_{k=1}^{n-1}\sum_{\sigma 
\in Sh(k,n-k)} (-1)^\sigma \epsilon(\sigma) v_{\sigma(1)} \wedge \cdots 
\wedge v_{\sigma(k)} \otimes v_{\sigma(k+1)} \wedge \cdots 
\wedge v_{\sigma(n)}, \nonumber
\end{eqnarray}
where $Sh(k,n-k)$ are  the 
unshuffles of type $(k,n-k)$, 
that is those permutations $\sigma$ of $n$ elements  that 
$\sigma(i) < \sigma(i+1) $ when $i \not= k.$ The sign 
$\epsilon(\sigma)$ is determined  by the requirement 
that 
\[
 v_1 \wedge \cdots 
\wedge v_n = (-1)^\sigma  \epsilon(\sigma) v_{\sigma(1)} \wedge \cdots 
\wedge v_{\sigma(n)}, 
\]
where $(-1)^\sigma$ is the  sign of the permutation $\sigma$.
Consider the suspension of the space $V; \ sV= V[1].$
\begin{dfn}
An $L_\infty$--algebra structure on a graded vector space $V$
is a codifferential $Q$ on $\Lambda(sV)$ of degree $+1,$ 
that is a map 
$Q: \Lambda (sV) \to \Lambda (sV)[1] $ such that
\begin{itemize}
\item $Q$ is a coderivation: 
$\Delta \circ  
Q  =  (Q \otimes 1 + 1 \otimes Q) \circ \Delta, $
\item  $Q \circ Q = 0.$
\end{itemize}
\end{dfn} 

 A coderivation $Q_k$ is of $k--$th order  if 
it is defined by a map $Q_k: \Lambda^k (sV) \to sV.$
Then the coderivation property provides the extension of the 
action of $Q_k$ on $\Lambda^n (sV)$ for any $n$:
\[
Q_k: \Lambda^n (sV) \to \Lambda^{n - k + 1} (sV) \quad \mbox{for} \quad n \geq  k, \
\quad \mbox{and} \quad Q_k: \Lambda^n (sV) \to 0 \quad \mbox{otherwise.}
\]
 This way we can consider sums of  coderivations of 
various orders  and define 
\begin{eqnarray*}
Q(v_1 \wedge \cdots 
\wedge v_n)&  \ \\ =
\sum_{k=1}^{n}  & \sum\limits_{\sigma \in Sh(k,n-k)}^{\mbox{ \  }} 
(-1)^\sigma \epsilon(\sigma) \  Q_k \bigl( v_{\sigma(1)} \wedge \cdots 
\wedge v_{\sigma(k)} \bigr) \wedge v_{\sigma(k+1)} \wedge \cdots 
\wedge v_{\sigma(n)},
\end{eqnarray*}
where $Q_k:  \Lambda^k (sV) \to sV$ and $
Q = \sum_{k=1}^{\infty} Q_k.$
 Then we can rewrite $Q^2 = 0$ as a sequence of equations for each $n$:
\[
\sum_{k=1}^{n} (-1)^{k(n-k)} \sum_{\sigma 
\in Sh(k,n-k)} (-1)^\sigma \epsilon(\sigma) Q_{n-k +1} 
\Bigl( Q_k \bigl(v_{\sigma(1)} \wedge \cdots 
\wedge v_{\sigma(k)}\bigr) \wedge v_{\sigma(k+1)} \wedge \cdots 
\wedge v_{\sigma(n)} \Bigr) = 0. 
\]

\begin{remark} 
\rm An  $L_\infty$--algebra $V$ has the following 
geometrical meaning.  For each $k:$
$\Lambda^k (sV)  = Sym^k V,$ $k-$th symmetric power of the 
space $V$. If $V$ is finite--dimensional, the 
symmetric powers of the space $ V$ 
are algebraic functions on the dual space $ V^\ast,$
which suggests that $Q$ be a vector field on the dual space.
$Q_k$ then  are   Taylor coefficients of the odd vector field $Q.$ 
 Hence the  map $Q$ could be interpreted 
as an odd vector field  of square $0.$ 
Such $Q$ is called a homological vector field.
 The notion of a homological vector field appears in \cite{V}, 
in relation to  the  Gerstenhaber structure on 
the exterior algebra of an  algebroid. 
A.S.Schwarz \cite{Schw} calls supermanifolds with a homological 
vector field $Q$--manifolds.
\end{remark}

\section{Deformations of Batalin--Vilkovisky algebras.}

The brackets (\ref{eq:FD})  are skew\-symmetric 
when the product $\mu$ is graded commutative. Hence they 
can be restricted to the exterior powers of $A:$
\[
F_D^k: \Lambda^kA \to A.
\]
We now  extend each linear map $F_D^k$ to a coderivation of $\Lambda A.$
We are going to show that the sum of all these coderivations 
is of square zero.

We need just another notion related to the degree:
\begin{dfn}
A linear map $D:A \to A$, where
$A = \sum_i A_i,$ is a $\Z$--graded vector space, is called {\it odd} if 
$D: A_i \to \sum_k A_{i+ 2k +1}, \ k \in \Z$ for each $i$.
\end{dfn}

\begin{proposition}\footnote{ While finishing this article, 
I learned about 
the paper \cite{BDA}
which contains a result similar to this proposition. 
However,  the aim and the language 
of \cite{BDA} are somewhat different.}
\label{th1}
Consider an odd operator $D$ on a 
graded commutative algebra $(A, \mu)$.
 Then $D^2 = 0$ if and only if 
the sum of brackets $Q_D = \sum F_D^n$ 
is a codifferential on $\Lambda A$ defining  an $ L_\infty$--structure,
in other words 
$ \sum_{k+l=n+1} F_D^k \circ F_D^l = 0$  for each $n \geq 1$.
\end{proposition}
\Proof

 The "if"  direction is obvious --- it is given by the 
first equation in the series of equations above: $n=k=l=1.$
The proof of the "only if" part is a tedious calculation.

For a graded commutative algebra, Akman's definition of the 
brackets (\ref{eq:FD})  coincides with the definition 
of Koszul \cite{Ko}, which we  reformulate 
in the following terms. 
 Define
a product on the exterior algebra
$M: A \wedge A \to A$ by $M (a_1 \wedge a_2) = a_1 \cdot a_2.$ 
We can extend it to any exterior power 
$ M (a_1 \wedge \ldots \wedge a_n) =
a_1 \cdot  \ldots \cdot a_n.$ Then we can define an 
$M$--coproduct as a map 
$\Lambda A \to A \otimes A: \ \Delta_M = (M \otimes M)\Delta:$
\[
\Delta_M (a_1 \wedge \ldots \wedge a_n) = 
\sum_{k=1}^{n-1}\sum_{\sigma 
\in Sh(k,n-k)} (-1)^\sigma \epsilon(\sigma) a_{\sigma(1)} \cdot \ldots 
\cdot  a_{\sigma(k)} \otimes a_{\sigma(k+1)} \cdot \ldots  \cdot
 a_{\sigma(n)}.
\]
Koszul's definition of multi-brackets is the following:
\[
 F_D^n  (a_1 \wedge \ldots \wedge a_n) = 
M (D\otimes 1) (a_1 \otimes 1 - 1 \otimes a_1)
\cdots (a_n \otimes 1 - 1 \otimes a_n).
\]
It can be reformulated as 
\begin{equation}
\label{eq:KFD}
F_D^n  (a_1 \wedge \ldots \wedge a_n) = M (D \otimes 1) \Delta_M 
(a_1 \wedge \ldots  \wedge a_n). 
\end{equation}
Then the lemma states that
\[
\Bigl( M (D \otimes 1) \Delta_M \Bigr) 
\Bigl(M (D \otimes 1) \Delta_M \otimes 1 \Bigr) \Delta = 0
\]
iff $D^2 =0.$
We see that, in the left hand side of this equation there are 
eather summands containing 
$D^2$ or summands which are present twice with opposite 
signs, due to the fact that the operator $D$ is odd. 
\endproof

Notice that the brackets  $F_D^n$ form an L$_\infty$ structure
with homotopies with respect to the operator $D,$ since 
the bracket $F_D^2$ gives a Lie algebra structure on $H(A,D),$
the cohomology of $A$ with respect to the operator $D.$

\begin{remark}
\label{rem:order}
\rm
{\it Order and degree.} 
There is a filtration on the algebra of differential operators
defined by their order.   For the operator $D$ however 
we would like to  obtain an unambiguous splitting 
$D= \Sigma _{n \geq 1} D_n,$ where $D_n$ are 
homogeneous operators of $n$-th order. 
All we know is  that for  the first $D_1,$ 
 $F^n_{D_1} \equiv 0 $ for $n > 1$. 
Then  $F^n_{D_2} \equiv 0, n>2,$ but  $F^2_{D_2} \not= 0,$ 
but there is already an ambiguity for ddefining $D_2.$

To obtain the splitting into homogeneous operators we use the 
the degree. 

$D$ acts on a graded algebra, so $D$ is a sum of operators of different 
degrees. It turns out that degree and order are in correspondence. 
It is natural to ask that 
the classical BV structure is a particular case of 
the generalized structure. Hence, we 
may start with the requirement that 
$D_1$ is of order $1$ and of
degree $+1$, and  $D_2$ is of order $2$ and  of degree $-1$. 
This defines the grading: the operator $D$
is unambiguously represented 
as a sum of  homogeneous 
operators. 
\begin{lemma}
\label{deg}
Consider an operator $D: A \to A,$ 
such that  $D^2 = 0$ and assume that $D$ is  the 
sum of  an operator 
 of order $1$ and of degree $+1, \ D_1: A_\bullet  \to A_{\bullet + 1}$ 
and  higher order  operators.
Then $D$ 
can be represented as a sum
 \[
D = \sum_{n \geq 1} D_n
 \]
where  each 
$D_n$ is an operator of  order $n$  and  of degree $3-2n,$ 
(in other words: $F^{n+1}_{D_n} \equiv 0$
and $D_n: A_\bullet  \to A_{\bullet +3-2n}.$) 
\end{lemma}
This lemma is an easy consequence of the condition $D^2 =0.$ 
 Of course we  can also weigh each operator of a certain degree
by some corresponding power of $\hbar.$
\end{remark}
\begin{remark}
\rm {\it Differential  BV--algebra.}
If the operator $D$ is of order $n$ we see that the highest homotopy
is given by  the $n$-th bracket. 

In particular, the second bracket 
\[
F_D^2 (a,b) = D (ab)  - Da~b - (-1)^{|a|} a~Db
\]
gives a classical BV--bracket for the case when $D_n = 0 $ for
$n \geq 3$.  Then the operator $D$ is of order $2$, that is
$D = D_1 + D_2.$ Such a $D$
describes the case of a differential BV--algebra which is 
the starting point of \cite{BK}, see also \cite{M}.

On the other hand, given  a differential 
algebra $(A, \mu, d)$  with additional second order 
differential operator $\delta$
one can define a generalized BV--algebra  by adding 
operators of higher order to $d+\delta$ , 
 requiring that their sum 
 \[
D = D_1 + D_2 + D_3 + \ldots
 \]
be of square $0,$ (here $D_1 =d, D_2 = \delta$).
Comparing with the differential BV--algebra case we see that 
there are still two differentials on the 
generalized algebra, $D$
and $D_1$ (the fact that $D_1$ is a differential  
 follows from $D^2 =0$).
 The following lemma is easy to prove.
\begin{lemma}
An operator on the algebra $(A,\mu), \ D = \sum D_n,$ such 
that $D^2 = 0$ 
is a derivation of the bracket 
$[a,b] = (-1)^{|a|} F_D^2 (a,b)$,
but not of  the product $\mu,$ 
while $D_1$ is a derivation of the product 
but not of the bracket.
\end{lemma}
\end{remark}

\begin{remark}\label{Leibniz}
\rm
{\it Generalization to Leibniz algebras}.
  If we start with a non-commu\-tative associative algebra structure,
the brackets $F^n_D$ (\ref{eq:FD}) still make sense for a
differential operator $D,$
(Definition \ref{Akman}).  However since there is no antisymmetry
condition anymore, the homotopy structure we get from
$D^2=0$ is not $L_\infty$.  Instead, one gets
Leib$_\infty$--algebra (\cite{Li}), 
homotopy version of a Leibniz
algebra  (\cite{L}). 
\end{remark}
\section{Commutative BV$_\infty$--algebra.}
We  now propose a definition of 
a strong homotopy Batalin--Vilkovisky algebra (BV$_\infty$--algebra).
 Here we will restrict ourselves to the case 
of commutative algebras. 
\begin{dfn}
\label{def}
A triple $(A,d,D)$ is a commutative BV$_\infty$--algebra when
\begin{itemize}
\item
$A$  is a graded commutative  algebra,
\item 
$d: A \to A$ is a degree $1$ differential of the algebra $A,$ 
\item
$D: A \to A$ is  an odd square zero differential operator, 
such that the degree 
of $D - d$ is negative.
\end{itemize}
\end{dfn}

There are various ways to define a
 BV$_\infty$--algebra. 
In our definition the commutative structure is preserved.
One can imagine deforming the commutative structure as well. In 
Remark \ref{Leibniz} we mentioned one of the generalizations, 
the one leading to the Leib$_\infty$--algebras.
However, all definitions should 
lead to the following property --- 
a  BV$_\infty$--algebra
should have a BV--algebra structure on its cohomology.
Indeed in our case:
\begin{theorem}
\label{th:2} 
The cohomology $H(A,d)$ of a commutative 
BV$_\infty$--algebra $(A,d,D)$  
is a BV--algebra.
\end{theorem}

\Proof
Consider the condition $D^2 = 0.$ Since the degree 
of $D - d $ is negative, it means that $D$ is the  sum of
$d,$ a derivation of degree $+1,$ and of negative degree operators:
$D_2 + D_3 + \cdots.$
From the fact that $D^2 = 0 $ follows that $d D_2 +D_2 d = 0,$ that is 
$D_2$ acts on the cohomology $H(A,d).$ Moreover, 
$D_2^2 = d D_3 +D_3 d,$ which means that on the cohomology 
$H(A,d), \ D_2^2 = 0.$  
Since  $D_2$ is a second order operator, it defines  
the structure of a  BV--algebra 
on the cohomology $H(A,d).$
\endproof

\section{Possible applications.}
 It would be interesting if we could 
extend the formality theorem of  Kontsevich
\cite{K} to the quasi-isomorphism of BV$_\infty$--algebras.

The formality theorem of Kontsevich claims that 
two differential graded Lie algebras
defined on any manifold $M,$ the  algebra of local 
Hochschild cochains and the 
algebra of polyvector fields, are quasi-isomorphic as 
$L_\infty$--algebras.

 Let  $A$ denote   the algebra of smooth functions on $M,$
 $A = C^\infty(M)$,
 with the  pointwise commutative product.
Let $D$ be the algebra of  polydifferential operators on 
$M\!: D = \oplus D^k, \ D^k  = Hom_{\mbox{\it loc}}(A^{\otimes k+1}, A),$
and  let $T$ be the algebra of polyvector fields on $M\!: T = \oplus T^k, 
\ T^k = \Gamma (\Lambda^{k+1}TM)$, both with the degree shifted by $1.$

\noindent Then there are the following corresponding structures on these two 
algebras: 
\vspace{0.5cm}

 \begin{center}
{\renewcommand{\arraystretch}{1,5}%
\begin{tabular}{||l|r|r||}
\hline
Graded space &\bf Polyvector fields &
\bf Polydifferential Operators \\
 &$ T = \oplus T^\bullet = \oplus  \Gamma(\Lambda^{\bullet+1} TM)$ & 
$D = \oplus D^\bullet = \oplus Hom_{\mbox{\it loc}} (A^{\bullet+1},A)$\\
\hline
 Differential & $ d=0$& Hochschild $B: D^\bullet \to D^{\bullet+1}$ \\
\hline
Lie bracket &   Schouten--Nijenhuis & Gerstenhaber\\ 
\hline
\strut Product &  $\wedge$ --- exterior product& $\cup$ --- cup product \\
\hline
\strut BV--operator &  $\delta$ & ?? \\
\hline 
\end{tabular}%
}
\end{center}

\vspace{0.5cm}

\noindent One can check that  $T$ is  
in fact a Gerstenhaber algebra  while 
$D$ is a Gerstenhaber algebra up to homotopy, since 
the $\cup$-product on $D$ is commutative only up to homotopy.
However the Lie adjoint action on $D$ 
is still an odd derivation with respect to the product. 

 Recently Dima Tamarkin \cite{T}  proved 
a generalization of Kontsevich's  formality  
theorem, he showed the existence of 
a morphism of Gerstenhaber algebras up 
to homotopy between $T$ and $D$. 
 In other words, the 
algebra of polydifferential operators is G--formal:
the algebra of polydifferential operators and 
 the algebra of polyvector fields are quasi-isomorphic as 
 G$_\infty$--algebras  (Gerstenhaber algebras up to homotopy).

We would like to see if one could prove the formality 
not only as  G$_\infty$--algebras but as  BV$_\infty$--algebras. 
 
If the first Chern class of a manifold $M$ is $0$, then
the algebra of polyvector fields on $M$ is  a BV--algebra.
There is a one-to-one correspondence between BV--structures 
on a manifold $M$ and 
flat connections on the determinant  bundle (bundle 
of polyvector fields in 
the top degree: $\Lambda^{{top}}TM$). 
Such a structure on real manifolds was studied in many papers
\cite{Ko,Xu,Hu,W}, on 
Calabi--Yau manifolds one can refer to \cite{Sch,BK}.
We conjecture that in these cases 
there should be some BV$_\infty$--structure
leading to the Gerstenhaber bracket on polydifferential
operators.  
\begin{conjecture}
There is a structure of a BV$_\infty$--algebra
on the space of polydifferential 
operators on a manifold with a 
zero first Chern class.
\end{conjecture} 
\begin{conjecture} The BV$_\infty$--algebra of  
polydifferential operators on a manifold is formal: it 
is quasi-isomorphic as  a BV$_\infty$--algebra to its cohomology, 
the BV--algebra of polyvector fields.
\end{conjecture} 

For these conjectures we will need a more general definition than
definition \ref{def}, since the cup product on the algebra of
polydifferential operators is commutative only up to homotopy.  This
generalization should not pose a problem, it will be done in a 
subsequent article.

From the conjecture, would follow the  
Maurer--Cartan equation (MC--equations) 
for the BV operator on the  algebra of polydifferential operators 
(probably tensored with some graded commutative algebra). Moreover,
a quasi-isomorphism of BV$_\infty$--algebras would map 
solutions of the  MC--equation on one algebra
to solutions of the MC--equation on the other  algebra.

We know from \cite{BK} that formal moduli space of 
solutions to the MC--equation, modulo gauge invariance
on polyvector fields tensored with the algebra of anti-holomorphic forms
on a Calabi--Yau manifold 
carries a natural structure of Frobenius manifold. 
If a quasi-isomorphism $ T \to D$ of BV--structures 
up to homotopy exists 
it would define a Frobenius manifold structure on the 
solutions of MC--equation modulo gauge invariance 
on polydifferential operators
tensored with the algebra of anti-holomorphic forms.

Another instance  where we could expect to find
generalized BV--structures is in the theory of vertex operator algebras.
There is a structure of a Batalin--Vilkovisky algebra on the
cohomology of  a vertex operator algebra
 (see  \cite{LZ}, \cite{PS}). 
It is natural to ask what  structure exists on the
vertex operator algebra itself. This shows the need  
for a suitable definition of a  BV$_\infty$--structure.
Besides it should fit into the 
general picture outlined by Stasheff \cite{S}.
 
{\bf Acknowledgments.} I am grateful to Paul Bressler, 
Ale Frabetti, Ezra Getzler, Johannes Huebschmann, Maxim Kontsevich,  
Yvette Kosmann-Schwarzbach, Jean-Louis Loday,
Dominique Manchon, 
A.S. Schwarz, Boris Tsygan, Sasha Voronov  and Simon Willerton
for inspiring discussions
and to the organizers of the Warsaw conference for their hospitality.  
 The style of this exposition owes a lot to the remarks of the referee.

\end{document}